\documentclass[reqno]{amsart}

\usepackage{graphicx,enumerate,booktabs,nicefrac,color,bm}
\usepackage{stmaryrd}
\usepackage{amsmath, amssymb, amstext, amsfonts, mathrsfs}
\usepackage{algorithm,algorithmicx,algpseudocode}
\usepackage{cite}

\newcommand{\norm}[1]{\left\|#1\right\|}

\newcommand{\abs}[1]{\lvert #1 \rvert}
\newcommand{\F}{\mathsf{F}}

\renewcommand{\d}{\mathsf{d}}

\newcommand{\dt}{\,\d t}
\newcommand{\dx}{\,\d\bm x}
\newcommand{\I}{\mathsf{I}_h^{n}}

\newcommand{\NN}[1]{\left|\!\left|\!\left|#1\right|\!\right|\!\right|} 

\newcommand{\solve}{\mathtt{solve}}
\newcommand{\dprod}[2]{\left<#1,#2\right>}

\newcommand{\myState}[1]{\State\parbox[t]{\dimexpr\linewidth-\algorithmicindent}{#1\strut}}

\newcommand{\myStateTriple}[1]{\State\parbox[t]{\dimexpr\linewidth-\algorithmicindent-\algorithmicindent-\algorithmicindent}{#1\strut}}

\newtheorem{Remark}[equation]{Remark}
\newenvironment{remark}{\begin{Remark}\rm}{\end{Remark}}
\newtheorem{theorem}[equation]{Theorem}

\newtheorem{lemma}[equation]{Lemma}
\newtheorem{ex}[equation]{Example}
\newenvironment{example}{\begin{ex}\rm}{\end{ex}}

\numberwithin{equation}{section}

\numberwithin{equation}{section}

\title[Adaptive Newton-Galerkin Time-Stepping Methods]{Fully Adaptive Newton-Galerkin Time Stepping Methods for Singularly Perturbed Parabolic Evolution Equations}

\author[M.~Amrein]{Mario Amrein}
\address{Mathematics Institute, University of Bern, CH-3012 Switzerland}
\email{mario.amrein@hslu.ch}
  
\author[T.~P.~Wihler]{Thomas P.~Wihler}
\address{Mathematics Institute, University of Bern, Sidlerstrasse 5, CH-3012 Switzerland}
\email{wihler@math.unibe.ch}

\begin{document}

\begin{abstract}
In this paper we develop an adaptive procedure for the numerical solution of semilinear parabolic problems, with possible singular perturbations. Our approach combines a linearization technique using Newton's method with an adaptive discretization---which is based on a spatial finite element method and the backward Euler time stepping scheme---of the resulting sequence of linear problems. Upon deriving a robust \emph{a posteriori} error analysis, we design a fully adaptive Newton-Galerkin time stepping algorithm. Numerical experiments underline the robustness and reliability of the proposed approach for various examples.
\end{abstract}



\maketitle

\section{Introduction}
Semilinear evolution problems appear in a wide range of applications including, for instance, ecology, (bio-)chemistry, quantum- and astro-physics, material science, or optics; see~\cite{BarlesBurdeau:95,BarlesPerthame:07,BorisyukErmentroutFriedmanTerman:05,CantrellCosner:03,Chandrasekhar:39,Edelstein-Keshet:05,JagerLuckhaus:92,Kelley:65,OkuboLevin:01,Perthame:07}. In this contribution we consider the numerical approximation of semilinear parabolic equations with a possibly degenerate diffusion coefficient~$0<\varepsilon\ll 1$. Specifically, for a continuously differentiable nonlinearity $ f:\mathbb{R}\times \Omega \times (0,T] \rightarrow \mathbb{R}$, and an initial function $g \in L^{\infty}(\Omega)$, we study initial/boundary value problems of the form
\begin{equation}
 \label{heat}
\left\{ \begin{aligned}
\partial_{t}u(\bm x,t) -  \varepsilon \Delta u(\bm x,t) &= f(u(\bm x,t),\bm x,t),  & &(\bm x,t)\in \Omega\times(0,T],\\
u(\bm x,t)&=0,  &&(\bm x,t)\in  \partial\Omega\times(0,T],\\
u(\bm x,0)&=g(\bm x),  & &\bm x\in\Omega.
\end{aligned} \right.
\end{equation}
Here $ \Omega \subset \mathbb{R}^{d} $, with $d\in \{1,2\} $, is an open and bounded one dimensional ($1$d) interval or a two-dimensional ($2$d) Lipschitz polygon, respectively. Furthermore, $T\in (0,\infty)$ denotes the final time of the evolutionary process. In the singularly perturbed case, solutions of~\eqref{heat} are known to exhibit boundary layers, interior shocks, or (multiple) spikes. The appearance of singular effects of this type constitutes a challenging issue when solving such problems numerically; see, e.g. \cite{RoStTo08,Verhulst}.

With the aim of designing an adaptive numerical procedure for~\eqref{heat}, we follow our recent work on stationary (elliptic) PDE in~\cite{AmreinWihler:14,AmreinWihler:15} (see also~\cite{Deuflhard:04,El-AlaouiErnVohralik:11}). In particular, this includes the application of a Newton linearization framework to the nonlinear problem at hand, and, subsequently, the discretization of the resulting \emph{sequence of linear evolution problems} by appropriate numerical schemes. It is worth mentioning that this methodology enables the use of numerical analysis techniques that originate from the treatment of~\emph{linear} problems; this is opposed to studying nonlinear discretizations schemes (see, e.g., \cite{georg,KyzaMakridakis:11}). The challenge in deriving practically effective \emph{a posteriori} error bounds within this setting is to provide a suitable splitting of the total residual into several computable quantities, each of which accounts for one of the different errors that have been committed during the discretization process: (1) a linearization residual, (2) a time discretization residual, and (3) a space discretization residual. Then, based on the resulting \emph{a posteriori} error estimates, a fully adaptive Newton-Galerkin time stepping algorithm for the numerical solution of~\eqref{heat} can be derived. Specifically, in order to obtain an efficient overall complexity of the scheme, we propose an interplay between the Newton linearization, time adaptivity, and spatial adaptivity. To do so, the algorithm will take into account the different residuals (1)--(3), and will perform a Newton step, a refinement of the current time step, or a refinement of the spatial mesh according to whatever residual is currently dominant. 

In the context of this paper, a $\mathbb{P}_1$-finite element approach in space, and a backward Euler discretization in time will be applied. Our \emph{a posteriori} error analysis proceeds along the lines of the theory presented in~\cite{Verfuerthbook} on linear parabolic equations. Furthermore, in order to obtain $\varepsilon$-robust bounds, we will follow the papers~\cite{Verfuerth} and~\cite{AmreinWihler:15} on finite element discretizations for singularly perturbed linear and semilinear elliptic problems, respectively. By means of a series of numerical experiments we will demonstrate that the interactive  application of temporal and spatial mesh refinements, together with a continued monitoring of the linearization effect, leads to an~$\varepsilon$-robust decay of the residual even in the singularly perturbed regime.

\subsection*{Outline}
The paper is organized as follows: In Section~\ref{sc:linearizatio/discretization} we begin by deriving the Newton linearization of~\eqref{heat}, and formulate the discretization of the resulting sequence of linearized problems in the spatial and temporal variables by means of a finite element method and the backward Euler scheme, respectively. Furthermore, the goal of Section~\ref{sec:aposteriori} is to derive an $\varepsilon$-robust {\em a posteriori} error analysis. Moreover, in Section~ \ref{sec:numerics}, we develop a fully adaptive Newton-Galerkin time stepping algorithm. Furthermore, we present a series of numerical experiments illustrating the performance of the proposed adaptive procedure. Finally, we summarize our findings in Section~\ref{sc:conclusions}.

\subsection*{Notation and Problem Formulation} 
For the purpose of this paper, we define the space~$Z:=H_{0}^{1}(\Omega)$, where $H^1_0(\Omega)$ is the standard Sobolev space of functions in $H^1(\Omega)~=W^{1,2}(\Omega)$, with zero trace on~$\partial\Omega$. The space~$Z$ is equipped with the singular perturbation norm~$\|\cdot\|_Z:=\NN{\cdot}_{\varepsilon,\Omega}$, where, for any subset~$D\subseteq\Omega$, we define
\[
\NN{u}_{\varepsilon,D}:=\Bigl(\varepsilon\norm{\nabla u}_{0,D}^2 +\norm{u}_{0,D}^2 \Bigr)^{\nicefrac{1}{2}},\qquad u\in H^1(D).
\]
Here, $\|\cdot\|_{0,D}$ denotes the $L^2$-norm on~$D$. Frequently, for~$D=\Omega$, the subindex~`$D$' will be omitted. In the sequel, we will abbreviate $f(u,\bm x,t)$ by $f(u)$; note that, in the case of~$f(u)=-u$, when~\eqref{heat} is linear, the norm $\NN{\cdot}_{\varepsilon,\Omega}$ is the natural energy norm on~$Z$.

Moreover, we signify by~$Z'=H^{-1}(\Omega)$ the dual space of~$Z$; it is equipped with the norm
\[
\|\psi\|_{Z'}=\sup_{\genfrac{}{}{0pt}{}{z\in Z}{\|z\|_Z=1}}{\dprod{\psi}{z}},
\]
where $ \dprod{\cdot}{\cdot}$ is the dual product in $Z'\times Z$. Furthermore, consider the Bochner spaces~$Y:=L^{2}(0,T;Z')$, and~$X:=\{u \in L^{2}(0,T;Z):\partial_t u \in Y\}$, with~$\partial_{t}$ being the time derivative operator in the distributional sense. On~$Y$ and $X$ we introduce the norms 
\[
\norm{u}_{Y}:=\left(\int_{0}^{T}{\norm{u(\cdot,t)}_{Z'}^{2}\dt}\right)^{\nicefrac{1}{2}},\qquad u\in Y,
\]
and
\[
\norm{u}_{X}:=\left(\int_{0}^{T}\left\{\|u\|_Z^2+\|\partial_tu\|_{Z'}^2\right\}\dt\right)^{\nicefrac{1}{2}},\qquad u\in X,
\]
respectively. 

Defining the map $\F_{\varepsilon}: X\rightarrow Y$ through
\begin{equation}
\label{eq:goal}
\dprod{\F_{\varepsilon}(u)}{v} :=  \dprod{\partial_{t}u}{v}+\int_{\Omega}\{\varepsilon \nabla u \cdot \nabla v-f(u)v\}\dx \qquad \forall v \in Z,
\end{equation}
the above problem \eqref{heat} can be written as a nonlinear operator equation in $Y$:
\[
 u \in X: \qquad \F_{\varepsilon}(u)=0 \qquad \forall t \in (0,T),
\]
with $ u(\cdot,0)=g $. 

Throughout this work we shall use the abbreviation $x \preccurlyeq y$ to mean $x \leq c y$, for a constant $c>0$ independent of the mesh size~$h$ and of~$\varepsilon>0 $. 


\section{Linearization and Discretization}\label{sc:linearizatio/discretization}
\subsection{Linearization}
We note that the Fr\'echet derivative of $ \F_{\varepsilon} $ from \eqref{eq:goal} at $u\in X$ is given by
\[
\left \langle \F_{\varepsilon}'(u)w,v\right \rangle = \dprod{\partial_{t}w}{v}+\int_{\Omega}{\{\varepsilon \nabla w \cdot \nabla v-\partial_{u}f(u)wv\}\dx},\qquad  w\in X,v\in Z.
\]
Then, starting from an initial guess $ u_0 \in X $, Newton's method for~\eqref{heat} is an iterative procedure by which we find $u_{N+1} \in X $ from $u_N \in X$, for~$N=0,1,2\ldots$, such that there holds  
\begin{equation}
\label{eq:linearization2}
\F_{\varepsilon}'(u_{N})(u_{N+1}-u_{N})=-\F_{\varepsilon}(u_{N})
\end{equation}
in $ Y$. Upon defining the increment $\delta_N:=u_{N+1}-u_{N}\in X$, and recalling~\eqref{eq:goal}, we note that
\begin{equation}
\label{eq:equiv}
\begin{split}
\dprod{\partial_{t}(u_{N}+\delta_{N})}{v}&+\int_{\Omega}\left\{\varepsilon \nabla(u_{N}+\delta_{N})\cdot \nabla{v}\right\}\dx\\
&=\int_{\Omega}{\{f(u_{N})+\partial_{u}f(u_{N})\delta_{N}\}v\dx},
\end{split}
\end{equation}
for all~$v\in Z$.

\subsection{Finite Element Meshes and Spaces}
   
Let $ \mathcal{T}_{h}=\{K\}_{K\in\mathcal{T}_h}$ be a regular and shape-regular mesh partition of $\Omega $ into disjoint open simplices, i.e., any~$K\in\mathcal{T}_h$ is an affine and non-degenerate image of the (open) reference simplex~$\widehat K=\{\widehat x\in\mathbb{R}_+^d:\,\sum_{i=1}^d\widehat x_i<1\}$. By~$h_K=\mathrm{diam}(K)$ we signify the element diameter of~$K\in\mathcal{T}_h$, and by $h=\max_{K\in\mathcal{T}_h}h_K$ the mesh size of~$\mathcal{T}_h$. Furthermore, by $\mathcal{E}_{h}$ we denote the set of all interior mesh nodes for~$d=1$ and interior (open) edges for~$d=2$ in~$\mathcal{T}_{h}$. In addition, for~$K\in\mathcal{T}_h$, we let~$\mathcal{E}_h(K)=\{E\in\mathcal{E}_h:\,E\subset\partial K\}$. For~$E\in\mathcal{E}_h$, we let~$h_E$ be the mean of the lengths of the adjacent elements in 1d, and the length of~$E$ in~2d.

We consider the finite element space of continuous, piecewise linear functions on $\mathcal{T}_{h}$ with zero trace on~$\partial\Omega$, given by
\begin{equation*}
V_{0}^{h}:=\{\varphi\in H^1_0(\Omega):\,\varphi|_{K} \in \mathbb{P}_{1}(K) \, \forall K \in \mathcal{T}_{h}\},
\end{equation*}
where~$\mathbb{P}_1(K)$ is the standard space of all linear polynomial functions on~$K$. Moreover, for any function $  \varphi \in V^h_0 $ and a given edge $ E \in \mathcal{E}_{h}$ with $E=\mathcal{E}_h(K^\sharp)\cap\mathcal{E}_h(K^\flat) $ shared by two neighboring simplices~$K^\sharp, K^\flat\in\mathcal{T}_h$, we denote by $\llbracket \varphi \rrbracket_{E} $ the (vector-valued) jump of $ \varphi $ across $E$:
\[
\llbracket \varphi \rrbracket_{E}(x)=\lim_{t\to 0^+}\varphi(\bm x+t\bm n^\sharp)\bm n^\sharp+\lim_{t\to 0^+}\varphi(\bm x+t\bm n^\flat)\bm n^\flat \qquad \forall \bm x \in E.
\]
Here, $\bm n^\sharp$ and~$\bm n^\flat$ denote the unit outward vectors on~$\partial K^\sharp$ and~$\partial K^\flat$, respectively. 

Furthermore, for any~$K\in\mathcal{T}_h$, 
we consider the element patch
\[
 \widetilde{w}_{K}:=\bigcup_{\genfrac{}{}{0pt}{}{K'\in\mathcal{T}_h:}{\overline K \cap \overline K'  \neq \emptyset}}{K'}.
\]
Let us also define the following two quantities:
\begin{equation}\label{boundary}
\begin{split}
\alpha_K&:=\min(1,\varepsilon^{-\nicefrac12}h_K),\qquad K\in\mathcal{T}_h,\\
\alpha_E&:=\min(1,\varepsilon^{-\nicefrac12}h_E),\qquad E\in\mathcal{E}_h.
\end{split}
\end{equation}
Then, we recall the following approximation result from~\cite[Corollary~4.2]{AmreinWihler:15}:

\begin{lemma}
Let $\I:\,H_{0}^{1}(\Omega)\rightarrow V_{0}^{h} $ be the quasi-interpolation Cl\'ement operator (see, e.g., \cite{Verfuerthbook}). Then, for any element~$K\in\mathcal{T}_h$, and any edge~$E\in\mathcal{E}_h$, with~$E=\mathcal{E}_h(K^\sharp)\cap\mathcal{E}_h(K^\flat)$ for some neighboring elements~$K^\sharp, K^\flat\in\mathcal{T}_h$, and any~$v\in H^1_0(\Omega)$, there hold the approximation bounds
\begin{align*}
\norm{v-\I v}_{0,K}&\preccurlyeq \alpha_{K}\NN{v}_{\varepsilon,\widetilde{w}_K},\\
\norm{v-\I v}_{0,E}&\preccurlyeq \frac12\varepsilon^{-\nicefrac14}\alpha_E^{\nicefrac12}\left(\NN{v}_{\varepsilon,\widetilde{w}_{K^\sharp}}+\NN{v}_{\varepsilon,\widetilde{w}_{K^\flat}}\right),
\end{align*}
where~$\alpha_K$ and~$\alpha_E$ are defined in~\eqref{boundary}.
\end{lemma}

\subsection{Newton-Galerkin Backward Euler Discretization}

In order to provide a numerical approximation of~\eqref{heat}, we will discretize the spatial and temporal variables in the formulation~\eqref{eq:linearization2} by means of a finite element method in space and the backward Euler scheme in time, respectively. In combination with the Newton iteration this results in a Newton-Galerkin time stepping approximation scheme. 

We consider a time partition of the interval $(0,T)$ into~$M\ge 1$ subintervals $ I_{n}=(t_{n-1},t_{n})$, $n=1,\ldots,M$, satisfying 
\[
0=t_{0}<t_{1}<\ldots<t_{M-1}<t_{M}=T,
\]
and define the time step lengths~$k_{n}:=t_{n}-t_{n-1}$. We mark any quantities related to the finite element discretization on a given time interval~$I_n$ by an index~`$n$'; in particular, we denote by $ \mathcal{T}_{h}^{n} $ the corresponding spatial partition of $\Omega$, and by 
\[
V_{0}^{h,n}:=\{\varphi\in H^1_0(\Omega):\,\varphi|_{K} \in \mathbb{P}_{1}(K) \, \forall K \in \mathcal{T}_{h}^{n}\}
\] 
the associated finite element space on a time subinterval~$I_n$. Furthermore, by $ \Pi^{n} $ we signify the $L^{2}$-projection onto $ V_{0}^{h,n} $. 

Applying the backward Euler time stepping scheme, the finite element discretization of the Newton iteration~\eqref{eq:equiv} on each time interval~$I_{n}$, $n=1,2,\ldots,M$, is to find $ \delta_{N}^{n} \in V_{0}^{h,n} $ from $ u_{N}^{n} \in V_{0}^{h,n} $ such that
\begin{equation}
\label{disc}
\begin{aligned}
\int_{\Omega}\bigg\{v\frac{u_{N}^{n}+\delta_{N}^{n}-u^{n-1}}{k_{n}}&+\varepsilon \nabla(u_{N}^{n}+\delta_{N}^{n})\cdot \nabla v\bigg\}\dx\\
=& \int_{\Omega}\left\{f^{n}(u_{N}^{n})+\partial_{u}f^{n}(u_{N}^{n})\delta_{N}^{n}\right\}v\dx \qquad \forall v \in V_{0}^{h,n},
\end{aligned}
\end{equation}
with the update $ u_{N+1}^{n}=u_{N}^{n}+\delta_{N}^{n}$, and with 
\[
f^{n}(\cdot):=f(u,\bm x,t)|_{(\cdot,\bm x,t_{n})}, \qquad 
\partial_{u}f^{n}(\cdot):=\partial_{u}f(u,\bm x,t)|_{(\cdot,\bm x,t_n)}.
\]
Moreover, for~$n=2,\ldots,M$, we denote by $u^{n-1}\in V_{0}^{h,n-1} $ the (space-dependent) discrete solution at the previous time node~$t_{n-1}$ (resulting from a sufficient number of Newton iterations), and for~$n=1$, we let $u^{0}:=\Pi^0g$ be a suitable approximation in~$H^1_0(\Omega)$ of the initial condition~$g\in L^\infty(\Omega)$. Furthermore, for the $n$-th time step, the initial guess~$u_{0}^{n}\in V_{0}^{h,n}$ is defined by
\begin{equation}\label{eq:u0}
u_{0}^{n}:=\Pi^{n}u^{n-1},\quad 1\leq n \leq M.
\end{equation}
We will denote the procedure of performing one Newton update, i.e., solving \eqref{disc} to obtain $ u_{N+1}^{n} $, by 
\[
u_{N+1}^{n}=\solve(k_{n},\mathcal{T}_{h}^{n},u_{N}^{n}).
\] 
Here, we make the assumption that we reinitiate the Newton iteration on each time step, i.e., for each $n=1,2,\ldots,M$, we start with $N=0$. 

\begin{remark}
\label{im}
The adaptive procedure based on the \emph{a posteriori} analysis to be presented in the next sections enables the use of possible coarsening of some (spatial) elements with small error contributions. 
A coarsening strategy within the procedure of solving \eqref{disc} is not a trivial task. In fact, suppose we have solved \eqref{disc} up to time $t_{n-1} $ so that we are given $u^{n-1}\in V_{0}^{h,n-1}$. Then, for a time step $k_{n}>0 $ small enough, it is reasonable to assume that $u^{n-1} $ 
is located in an attracting $\epsilon$-ball $ B_{\epsilon}(u_{\infty}^{n})\subset \mathcal{A}(u_{\infty}^{n}) $ of $ u_{\infty}^{n}\in V_{0}^{h,n} $, where $\mathcal{A}(u_{\infty}^{n}) $ denotes the attractor of the Newton iteration corresponding to  $u_{\infty}^{n}$.
Thence, we have
\begin{equation}
\label{attra}
\norm{u_{0}^{n}-u_{\infty}^{n}}
\leq \norm{\Pi^{n}u^{n-1}-u^{n-1}} + \norm{u^{n-1}-u_{\infty}^{n}} \leq \norm{\Pi^{n}u^{n-1}-u^{n-1}} + \epsilon,
\end{equation}
for some suitable norm $\norm{\cdot} $.
Hence, if $\mathcal{T}_{h}^{n}$ is obtained from $\mathcal{T}_{h}^{n-1}$ by refinement only, then  we see from \eqref{attra} that $ u_{0}^{n}\in \mathcal{A}(u_{\infty}^{n}) $ since $\Pi^{n}u^{n-1}=u^{n-1}$. If, however, 
there is a partial coarsening involved, we usually have 
$\Pi^{n}u^{n-1}\neq u^{n-1} $, and in consequence, the quantity $\norm{\Pi^{n}u^{n-1}-u^{n-1}}$ in \eqref{attra} may be too large in order to guarantee for
$u_{0}^{n} $ to stay within $ \mathcal{A}(u_{\infty}^{n}) $. Thus any coarsening strategy should only remove those degrees of freedom for which
$ \norm{\Pi^{n}u^{n-1}-u^{n-1}} $ remains of moderate size (cf.~also \cite{georg}).
\end{remark}

\section{\emph{A Posteriori} Error Analysis}\label{sec:aposteriori}

The goal of this section is to derive a residual based \emph{a posteriori} error bound for the discretization scheme~\eqref{disc}, which can be employed for the purpose of formulating an adaptive refinement procedure for the meshes and time steps in each Newton step. This leads to a fully adaptive Newton-Galerkin backward Euler discretization method for~\eqref{heat}. In the subsequent \emph{a posteriori} error analysis we follow closely the approach presented in \cite{Verfuerthbook}.

\subsection{Residuals}

The discrete problem \eqref{disc} generates a sequence $\{u_{N}^{n}\}_{N\ge 0}$ for each time step~$n=1,\ldots,M$. We can thus define a function~$u_{\mathcal{I}}\in C^0([0,T];H_{0}^{1}(\Omega))$ by
\[
u_{\mathcal{I}}|_{I_{n}}:=\frac{t_n-t}{k_n}u^{n-1}+\frac{t-t_{n-1}}{k_{n}}u_{N+1}^{n}, \qquad t \in [t_{n-1},t_{n}].
\]
We remark that $u_{\mathcal{I}}$ is understood as a function in time that depends on the (varying) Newton iteration index~$N$ on each subinterval~$I_n$. For later purposes notice that 
\begin{equation}
\label{brevity}
 q_{n}(t)\left(u^{n-1}-u_{N+1}^{n}\right)=u_{\mathcal{I}}-u_{N+1}^{n},\end{equation}
where $q_{n}(t):=k_n^{-1}(t_n-t)$, and, moreover, we observe that 
\[
\partial_{t} u_{\mathcal{I}}|_{I_{n}}=\frac{u_{N+1}^{n}-u^{n-1}}{k_{n}},
\]
for~$n=1,\ldots,M$. Therefore, motivated by the linear case discussed in \cite{Verfuerthbook}, we decompose the residual $ \F_{\varepsilon}(u_{\mathcal{I}}|_{I_n}) $ from~\eqref{eq:goal} on each time interval $I_{n}$, $n=1,\ldots,M$, as
\begin{equation}
\label{decomposition}
\left \langle \F_{\varepsilon}(u_{\mathcal{I}}|_{I_n}),v \right \rangle= \left \langle \F_{\varepsilon}^{1}(u_{\mathcal{I}}|_{I_n}),v \right \rangle+ 
\left \langle \F_{\varepsilon}^{2}(u_\mathcal{I}|_{I_n}),v \right \rangle+\left \langle \F_{\varepsilon}^{3}(u_\mathcal{I}|_{I_n}),v \right \rangle,\qquad v\in Z.
\end{equation}
Here, the parts $\F_{\varepsilon}^{i}(u_{\mathcal{I}}|_{I_n})$, $i =1,2,3$, are defined by
\begin{align*}
\left \langle \F_{\varepsilon}^{1}(u_{\mathcal{I}}|_{I_n}),v \right \rangle&:=\int_{\Omega}{\{v\partial_{t} u_{\mathcal{I}}|_{I_n}+\varepsilon \nabla u_{N+1}^{n}\cdot\nabla v-(f^n(u_{N}^{n})+\partial_{u}f^{n}(u_{N}^{n})\delta_{N}^{n})v\}\dx},\\
\left \langle \F_{\varepsilon}^{2}(u_{\mathcal{I}}|_{I_n}),v \right \rangle&:=\int_{\Omega}\varepsilon \nabla(u_{\mathcal{I}}|_{I_n}-u_{N+1}^{n})\cdot\nabla v\dx+  \int_{\Omega}{\{f^{n}(u_{N+1}^{n})-f(u_{\mathcal{I}}|_{I_n})\}v\dx},\\
\left \langle \F_{\varepsilon}^{3}(u_{\mathcal{I}}|_{I_n}),v \right \rangle&:=\int_{\Omega}{\{f^{n}(u_{N}^{n})+\partial_{u}f^{n}(u_{N}^{n})\delta_{N}^{n}-f^{n}(u_{N+1}^{n})\}v\dx},
\end{align*}
for any $v \in Z$. We emphasize that this splitting is crucial when aiming at an adaptive algorithm that is able to identify the individual error contributions resulting from the time and space discretizations as well as from the Newton linearization. In accordance with the notation introduced in \cite{Verfuerthbook}, we call $\F_{\varepsilon}^{1}(u_{\mathcal{I}}|_{I_n})$ and $\F_{\varepsilon}^{2}(u_{\mathcal{I}}|_{I_n}) $ the \emph{spatial} and the \emph{temporal residuals}, respectively. 
Furthermore, $\F_{\varepsilon}^{3}(u_{\mathcal{I}}|_{I_n})$ is termed the \emph{nonlinear residual}.

\subsection{\emph{A Posteriori} Error bound} 
Upon applying the triangle inequality to the decomposition \eqref{decomposition} we obtain that
\begin{equation}\label{eq:triangle}
\begin{split}
&\norm{\F_{\varepsilon}(u_{\mathcal{I}}|_{I_n})}_{L^{2}(I_{n};Z')}\\
&\quad\leq  \norm{\F_{\varepsilon}^{1}(u_{\mathcal{I}}|_{I_n})}_{L^{2}(I_{n};Z')}
+\norm{\F_{\varepsilon}^{2}(u_{\mathcal{I}}|_{I_n})}_{L^{2}(I_{n};Z')}+\norm{\F_{\varepsilon}^{3}(u_{\mathcal{I}}|_{I_n})}_{L^{2}(I_{n};Z')},
\end{split}
\end{equation}
on each time interval $I_{n}$, $n=1,\ldots,M$. We will now derive individual error bounds for each of the three residual terms $\F_{\varepsilon}^{i}(u_{\mathcal{I}}|_{I_n})$, $i=1,2,3$.
\subsubsection*{Spatial Residual}
We note the fact that the spatial residual~$\F_{\varepsilon}^{1}(u_{\mathcal{I}}|_{I_n})$ is constant with respect to time. It can thus be estimated as in the stationary case~\cite[Theorem~4.4]{AmreinWihler:15}. In fact, observing that
\[
\norm{\F_{\varepsilon}^{1}(u_{\mathcal{I}}|_{I_n})}_{L^{2}(I_n;Z')}=\sqrt{k_n}\|\F_{\varepsilon}^{1}(u_{\mathcal{I}}|_{I_n})\|_{Z'},
\]
we infer the estimate
\begin{equation}
\label{esf}
\norm{\F_{\varepsilon}^{1}(u_{\mathcal{I}}|_{I_n})}_{L^{2}(I_n;Z')}^2 \preccurlyeq k_{n}\sum_{K\in \mathcal{T}_{h}^{n}}{\eta_{n,K,N}^{2}},
\end{equation}
for any time interval $I_n$, $n=1,\ldots,M$. Here, for any $K\in \mathcal{T}_{h}^{n} $, the quantities
\begin{equation}\label{eq:eta}
\begin{aligned}
\eta_{n,K,N}^{2}:&=\alpha_{K}^{2}\norm{f^{n}(u_{N}^{n})+\partial_{u}f^{n}(u_{N}^{n})\delta_{N}^{n}+\varepsilon \Delta u_{N+1}^{n}-\partial_{t}u_{\mathcal{I}}|_{I_n}}_{0,K}^2\\
&\quad+\frac{1}{2}\sum_{E\in \mathcal{E}_{h}^{n}(K)}{\varepsilon^{-\nicefrac12}\alpha_{E}\norm{\varepsilon \llbracket \nabla u_{N+1}^{n}\rrbracket}_{0,E}^{2}}
\end{aligned}
\end{equation}
are computable error indicators, with~$\alpha_K$ and~$\alpha_E$ being defined in~\eqref{boundary}. We emphasize that the bound~\eqref{esf} is robust with respect to the singular perturbation parameter~$\varepsilon$.

\subsubsection*{Temporal Residual} 
Using the identity \eqref{brevity} we have that 
\[
\begin{aligned}
\left \langle \F_{\varepsilon}^{2}(u_{\mathcal{I}}|_{I_n}),v \right \rangle &= q_n(t)\int_{\Omega}{\{\varepsilon \nabla (u^{n-1}-u_{N+1}^{n})\cdot \nabla v\}\dx}\\
&\quad+\int_{\Omega}{\{f^{n}(u_{N+1}^{n})-f(u_{\mathcal{I}}|_{I_n})\}v\dx},
\end{aligned}
\]
from which, by application of the Cauchy-Schwarz inequality, we obtain
\[
\|\F_{\varepsilon}^{2}(u_{\mathcal{I}}|_{I_n})\|_{Z'}^{2} \leq q_{n}(t)^{2}\varepsilon\|\nabla(u^{n-1}-u_{N+1}^{n})\|_0^{2}+\norm{f^{n}(u_{N+1}^{n})-f(u_{\mathcal{I}}|_{I_n})}_{0}^{2}. 
\]
Moreover, since $\int_{t_{n-1}}^{t_n}{q_{n}(t)^2dt} = \nicefrac{k_n}{3}$,
we arrive at the bound
\begin{equation}
\label{eq:time}
\begin{split}
\norm{\F_{\varepsilon}^{2}(u_{\mathcal{I}}|_{I_n})}_{L^{2}(I_{n};Z')}^{2} &\leq  \frac{\varepsilon k_{n}}{3}\|\nabla(u^{n-1}-u_{N+1}^{n})\|_{0}^{2}\\
&\quad+\norm{f^{n}(u_{N+1}^{n})-f(u_{\mathcal{I}}|_{I_n})}_{L^{2}(I_{n};L^{2}(\Omega))}^{2},
\end{split}
\end{equation}
on each time interval $I_{n}$, $n=1,\ldots,M$. 

\begin{remark}
Applying the trapezoidal rule we see that we can approximate
\[
\norm{f^{n}(u_{N+1}^{n})-f(u_{\mathcal{I}}|_{I_n})}_{L^{2}(I_{n};L^{2}(\Omega))}^{2}
\approx k_{n}\norm{f^{n}(u_{N+1}^{n})-f^{n-1}(u^{n-1})}_{0}^{2},
\]
up to an error of order~$\mathcal{O}(k_{n}^2)$, for each time interval $I_{n}$, $n=1,\ldots,M$.
\end{remark}

\subsubsection*{Nonlinear Residual}
We immediately infer
\[
\|\F_{\varepsilon}^{3}(u_{\mathcal{I}}|_{I_n})\|_{Z'} \leq \norm{f^{n}(u_{N}^{n})+\partial_{u}f^{n}(u_{N}^{n})\delta_{N}^{n}-f^{n}(u_{N+1}^{n})}_0,
\]
and hence,
\begin{equation}
\label{nonlin}
\norm{\F_{\varepsilon}^{3}(u_{\mathcal{I}}|_{I_n})}_{L^{2}(I_{n};Z')}^2 \leq k_n\norm{f^{n}(u_{N}^{n})+\partial_{u}f^{n}(u_{N}^{n})\delta_{N}^{n}-f^{n}(u_{N+1}^{n})}_{0}^2,
\end{equation}
on each $I_{n}$, $n=1,\ldots,M$. 

Combining the bounds~\eqref{eq:triangle}, \eqref{esf}, \eqref{eq:time}, and~\eqref{nonlin} leads to the following result.

\begin{theorem}
\label{lam:1}
On each time interval $I_{n}$, $n=1,2,\ldots,M$, there holds the \emph{a posteriori} error bound
\[
\norm{\F_{\varepsilon}(u_{\mathcal{I}}|_{I_n})}_{L^{2}(I_n;Z')}^2
 \preccurlyeq   
  k_{n}\sum_{K\in \mathcal{T}_{h}^{{n}}}\left\{
  \eta_{n,K,N}^{2}+
  \vartheta^{2}_{n,K,N} + \Upsilon^{2}_{n,K,N}\right\},
\]
where, for~$K\in\mathcal{T}^n_h$, we recall the spatial error indicators~$\eta_{n,K,N}$ in~\eqref{eq:eta}, and let
\begin{equation}
\label{eq:definition}
\begin{aligned}
\vartheta_{n,K,N}^{2}&:= k_n^{-1}\norm{f^{n}(u_{N+1}^{n})-f(u_{\mathcal{I}}|_{I_n})}_{L^{2}(I_{n};L^{2}(K))}^{2}+\frac{1}{3}\varepsilon\|\nabla(u^{n-1}-u_{N+1}^{n})\|_{0,K}^{2},\\
\Upsilon_{n,K,N}^{2}&:=\norm{f^{n}(u_{N}^{n})+\partial_{u}f^{n}(u_{N}^{n})\delta_{N}^{n}-f^{n}(u_{N+1}^{n})}_{0,K}^2.
\end{aligned}
\end{equation}
\end{theorem}

For later reference, in addition to the previously introduced local \emph{a posteriori quantities}, we define the corresponding global error indicators as follows:
\begin{equation}\label{eq:global}
\begin{split}
\eta_{n,\Omega,N}^2&:=\sum_{K\in\mathcal{T}^n_h}\eta^2_{n,K,N},\qquad
\vartheta_{n,\Omega,N}^2:=\sum_{K\in\mathcal{T}^n_h}\vartheta^2_{n,K,N},\\
\Upsilon_{n,\Omega,N}^2&:=\sum_{K\in\mathcal{T}^n_h}\Upsilon^2_{n,K,N}.
\end{split}
\end{equation}

\subsection{Residual and Error Norm}\label{sc:RN}
Under certain conditions on the nonlinearity~$f$ in~\eqref{heat} it can be shown that the residual $\F_{\varepsilon}(u_{\mathcal{I}}) $ defined in \eqref{eq:goal} and the error $u- u_{\mathcal{I}}$ are equivalent. For example,
suppose that the nonlinearity~$f$ is Lipschitz continuous with Lipschitz constant~$L>0$, i.e.,
\[
|f(u,\bm x,t)-f(v,\bm x,t)|\le L|u-v|,
\]
and that it satisfies the monotonicity condition
\[
(f(u,\bm x,t)-f(v,\bm x,t))(u-v)\leq 0,
\]
for all $u,v\in\mathbb{R}, \bm x\in\Omega, t\in[0,T]$. Then, if $u$ is the exact solution of \eqref{heat}, we have that
\[
\begin{aligned}
-\langle \F_{\varepsilon}(u_{\mathcal{I}}),u-u_{\mathcal{I}}  \rangle 
&=  \frac{1}{2}\frac{\d}{\dt}\norm{u-u_{\mathcal{I}}}_{0}^{2}+\varepsilon \norm{\nabla(u-u_{\mathcal{I}})}_{0}^{2}\\
&\quad-\int_{\Omega}{(f(u)-f(u_{\mathcal{I}}))(u-u_{\mathcal{I}})\dx} \\
&\geq \frac{1}{2}\frac{\d}{\dt}\norm{u-u_{\mathcal{I}}}_{0}^{2}+\varepsilon \norm{\nabla(u-u_{\mathcal{I}})}_{0}^{2},
\end{aligned}
\]
for any~$t\in(0,T)$. Proceeding as in~\cite[Eq.~(4.13)]{AmreinWihler:15} it holds that
\[
-\left \langle \F_{\varepsilon}(u_{\mathcal{I}}),u-u_{\mathcal{I}} \right \rangle 
\geq \frac{1}{2}\frac{\d}{\dt}\norm{u-u_{\mathcal{I}}}_{0}^{2}+C_\varepsilon\NN{u-u_{\mathcal{I}}}_{\varepsilon}^{2},
\]
for a constant~$C_\varepsilon>0$. Furthermore, choosing $\delta >0 $ such that  $C_{\delta,\varepsilon}:=C_\varepsilon-\nicefrac{\delta}{2}>0$, and
using that
\[
\left|\left \langle \F_{\varepsilon}(u_{\mathcal{I}}),u-u_{\mathcal{I}} \right \rangle\right| \leq \frac{1}{2}\left(\delta^{-1}\|\F_{\varepsilon}(u_{\mathcal{I}})\|_{Z'}^{2}+\delta\NN{u-u_{\mathcal{I}}}_{\varepsilon}^{2} \right),
\]
we conclude
\begin{equation}
\label{imp}
\frac{1}{2}\frac{\d}{\dt}\norm{u-u_{\mathcal{I}}}_{0}^{2}+C_{\delta,\varepsilon}\NN{u-u_{\mathcal{I}}}_{\varepsilon}^{2}  
\leq \frac{1}{2\delta}\|\F_{\varepsilon}(u_{\mathcal{I}})\|_{Z'}^{2}.
\end{equation}
Integrating with respect to~$t$ yields
\begin{equation}
\label{final1}
\begin{aligned}
\int_0^t\NN{u-u_{\mathcal{I}}}_{\varepsilon}^{2}\dt
& \leq  \frac{\max{\{1,\delta^{-1}\}}}{2C_{\delta,\varepsilon}}\left(\norm{\F_{\varepsilon}(u_{\mathcal{I}})}_{L^{2}(0,t;Z')}^{2}+\norm{g-\Pi^{0}g}_{0}^{2}\right).
\end{aligned}
\end{equation}
In addition, invoking again \eqref{imp}, we obtain
\begin{equation}\label{eq:final2}
\norm{u-u_{\mathcal{I}}}_{L^{\infty}(0,t;L^{2}(\Omega))}^{2}
\leq \max\{1,\delta^{-1}\} \left(\norm{\F_{\varepsilon}(u_{\mathcal{I}})}_{L^{2}(0,t;Z')}^{2}+\norm{g-\Pi^{0}g}_{0}^{2}\right).
\end{equation}
Moreover, applying the Lipschitz continuity of~$f$, we observe, for $ v\in Z $ and~$t\in(0,T]$, that
\[
\abs{ \langle \partial_{t}(u-u_{\mathcal{I}}),v  \rangle }
 \leq  \abs{\left \langle \F_{\varepsilon}(u_{\mathcal{I}}),v \right \rangle }+\int_{\Omega}{\{\varepsilon \abs{\nabla(u-u_{\mathcal{I}})}\abs{\nabla{v}}
+L\abs{u-u_{\mathcal{I}}}\abs{v}\}\dx,}
\]
i.e., there holds
\[ 
\abs{\left \langle \partial_{t}(u-u_{\mathcal{I}}),v \right \rangle } \leq \abs{\left \langle \F_{\varepsilon}(u_{\mathcal{I}}),v \right \rangle}+\max{(1,L)}\NN{v}_{\varepsilon}\NN{u-u_{\mathcal{I}}}_{\varepsilon}.
\]
Thus,
\begin{equation}\label{eq:aux1}
\|\partial_{t}(u-u_{\mathcal{I}})\|_{Z'}\leq \|\F_{\varepsilon}(u_{\mathcal{I}})\|_{Z'}+\max{(1,L)}\NN{u-u_{\mathcal{I}}}_{\varepsilon}.
\end{equation}
Taking the square in the previous inequality, integrating over $(0,t)$, and recalling \eqref{final1} leads to
\begin{equation}
\label{her}
\norm{\partial_{t}(u-u_{\mathcal{I}})}_{L^{2}(0,t;Z')}^{2}\le C_{\delta,\varepsilon,L}\left(\norm{\F_{\varepsilon}(u_{\mathcal{I}})}_{L^{2}(0,t;Z')}^{2}+\norm{g-\Pi^{0}g}_{0}^{2}\right),
\end{equation}
for a constant~$C_{\delta,\varepsilon,L}>0$.
Combining \eqref{final1}, \eqref{eq:final2}, and \eqref{her}, we finally get that
\begin{align*}
E(t;u_{\mathcal{I}},g):&=\norm{g-\Pi^{0}g}_{0}^2+\norm{u-u_{\mathcal{I}}}_{L^{\infty}(0,t;L^{2}(\Omega))}^{2}\\
&\quad+\int_0^t\left\{\NN{u-u_{\mathcal{I}}}_{\varepsilon}^{2}+\|\partial_{t}(u-u_{\mathcal{I}})\|^2_{Z'}\right\}\dt\\
&\le \widetilde C_{\delta,\varepsilon,L} \left(\norm{\F_{\varepsilon}(u_{\mathcal{I}})}_{L^{2}(0,t;Z')}^{2}+\norm{g-\Pi^{0}g}_{0}^{2}\right),
\end{align*}
for a constant~$\widetilde C_{\delta,\varepsilon,L}>0$, and any~$t\in(0,T]$.
Conversely, proceeding as in~\eqref{eq:aux1}, it is possible to show that
\[
\|\F_{\varepsilon}(u_{\mathcal{I}})\|_{Z'}\leq \|\partial_{t}(u-u_{\mathcal{I}})\|_{Z'}+\max{(1,L)}\NN{u-u_{\mathcal{I}}}_{\varepsilon}.
\]
Integrating over~$(0,t)$ shows the equivalence of the residual term~$\norm{\F_{\varepsilon}(u_{\mathcal{I}})}_{L^{2}(0,t;Z')}+\norm{g-\Pi^{0}g}_{0}$ and of the error~$E(t;u_{\mathcal{I}},g)$.

\subsection{A Fully Adaptive Newton-Galerkin Algorithm}
We will now propose a procedure that will combine a Newton method with automatic spatial finite element mesh and time step refinements based on the \emph{a posteriori} error estimate from Theorem~\ref{lam:1}. Recalling our derivations in the previous Section~\ref{sc:RN}, it is reasonable to control the quantity
\[
E^n(u_{\mathcal I},g):=\norm{\F_{\varepsilon}(u_{\mathcal{I}})}^2_{L^{2}(0,t_n;Z')}+\norm{g-\Pi^{0}g}_{0}^2,
\]
for~$n=1,\ldots,M$. Then, by means of Theorem~\ref{lam:1}, we have that
\begin{equation}\label{eq:En}
E^n(u_{\mathcal I},g)
\preccurlyeq \eta_0^2+\sum_{j=1}^nk_j\left(\eta_{j,\Omega,N}^2+\vartheta^2_{j,\Omega,N}+\Upsilon^2_{j,\Omega,N}\right),
\end{equation}
where~$\eta_0:=\|g-\Pi^0g\|_0$, and~$\eta_{j,\Omega,N}$, $\vartheta_{j,\Omega,N}$, and~$\Upsilon_{j,\Omega,N}$ are defined in~\eqref{eq:global}. Given a final time $T>0 $, and some positive tolerances $ \varepsilon_{0}, \varepsilon_{\eta}, \varepsilon_{\vartheta}, \varepsilon_{\Upsilon}>0$, we define 
\begin{equation}\label{eq:localtolerances}
\varepsilon_{\text{loc},\eta}:=\frac{\varepsilon_{\eta}}{\sqrt{T}},\quad  \varepsilon_{\text{loc},\vartheta}:=\frac{\varepsilon_{\vartheta}}{\sqrt{T}},\quad \varepsilon_{\text{loc},\Upsilon}:=\frac{\varepsilon_{\Upsilon}}{\sqrt{T}}.
\end{equation}
Suppose that
\begin{equation}\label{eq:eta0}
\eta_0\le \varepsilon_0,
\end{equation} 
and that, for any~$n=1,\ldots,M$, there holds
\[
\eta_{n,\Omega,N}\le\varepsilon_{\text{loc},\eta},\quad
\vartheta_{n,\Omega,N}\le\varepsilon_{\text{loc},\vartheta},\quad
\Upsilon_{n,\Omega,N}\le\varepsilon_{\text{loc},\Upsilon}.
\]
Then, we conclude
\begin{equation}
\label{Finale}
E^M(u_{\mathcal I},g)\preccurlyeq \varepsilon_{T}^2,
\end{equation}
where
\[ 
\begin{aligned}
\varepsilon_{T}^{2}:=\varepsilon_{0}^{2}+\varepsilon_{\eta}^{2}+\varepsilon_{\vartheta}^{2}+\varepsilon_{\Upsilon}^{2}.
\end{aligned}
\]

We now formulate a possible realization of a time-space-Newton-Galerkin adaptive algorithm which aims to generate a numerical solution $ u_{\mathcal{I}} $ that satisfies the error bound~\eqref{Finale}. The basic idea is to exploit the structure of the error bound~\eqref{eq:En} to provide an interplay between adaptive finite element space refinements (and derefinements), automatic selection of the time steps, and an appropriate resolution of the Newton linearization error. More precisely, our adaptive procedure identifies whichever of the computable \emph{a posteriori} quantities occurring in~\eqref{eq:En} is currently dominant, and performs a corresponding refinement. In this way, the scheme follows along the lines of our previous approach in~\cite{AmreinWihler:14,AmreinWihler:15} on stationary problems, with the additional feature that the temporal errors are now taken into account, too. Our fully adaptive Newton-Galerkin method is outlined in Algorithm~\ref{spacetime}.

\begin{algorithm}
\caption{Fully-adaptive Newton-Galerkin time stepping method}
\label{spacetime}
\begin{algorithmic}[1]
\State Initialization: Input a final time $T>0$. Prescribe the overall error~$\varepsilon_{T}>0$, a lower bound for the time steps~$k_{\min}>0$, an initial time step~$k_{0}\ge k_{\min}$, and two time mesh parameters~$\kappa>1$ (derefinement) and $\sigma \in (0,1)$ (refinement). Set~$n\gets 1$.
\State Find an approximation~$\Pi^0g\in H^1_0(\Omega)$ of the initial condition~$g$ such that~$\varepsilon_0:=\|g-\Pi^0\|_0<\varepsilon_T$.
\State Split the tolerance~$\varepsilon_{T}$ into contributions~$\varepsilon_{\eta}, \varepsilon_{\vartheta}, \varepsilon_{\Upsilon}>0$ such that there holds $\varepsilon_{T}^{2}:=\varepsilon_{0}^{2}+\varepsilon_{\eta}^{2}+\varepsilon_{\vartheta}^{2}+\varepsilon_{\Upsilon}^{2}$.
\While {$k_{n}:=\min\{k_{n-1},T-t_{n-1}\}\ge k_{\min}$} \Comment{time stepping}
\myState {Generate a first mesh~$\mathcal{T}_{h}^{n} $ on the current, $n$-th, time step, and compute the initial guess~$u_0^n$ for the Newton iteration from~\eqref{eq:u0};  set $N\gets 0$.}
\myState {Compute $ u_{N+1}^{n}=\solve(k_{n},\mathcal{T}_{h}^{n},u_{N}^{n}) $, and evaluate the error indicators from~\eqref{eq:eta} and~\eqref{eq:definition}.}
\If {$\eta_{n,\Omega,N}^2+\vartheta_{n,\Omega,N}^2+\Upsilon_{n,\Omega,N}^2
> \varepsilon_{\text{loc},\eta}^{2}+\varepsilon_{\text{loc},\vartheta}^{2}+\varepsilon_{\text{loc},\Upsilon}^{2}$}
\If {$\vartheta_{n,\Omega,N}^2+\Upsilon^2_{n,\Omega,N}<\eta_{n,\Omega,N}^2$} 
\myStateTriple{refine the current mesh~$\mathcal{T}_{h}^{n}$ according to the elemental contributions to obtain a new mesh~$\mathcal{T}_{h}^{n}\gets\mathcal{T}_{h}^{n}$; go back to step (6). }
\ElsIf {$\vartheta_{n,\Omega,N}^2+\Upsilon^2_{n,\Omega,N}\ge\eta_{n,\Omega,N}^2$ and $\Upsilon_{n,\Omega,N}< \vartheta_{n,\Omega,N}$} 
\myStateTriple{set $ k_{n}\gets\sigma k_{n}$, $N\gets 0$, and go to step $(5)$ (if~$k_n\ge k_{\min}$, otherwise stop algorithm).}
\ElsIf {$\vartheta_{n,\Omega,N}^2+\Upsilon^2_{n,\Omega,N}\ge\eta_{n,\Omega,N}^2$ and $\Upsilon_{n,\Omega,N}\ge \vartheta_{n,\Omega,N}$} 
\myStateTriple {do another Newton iteration by going back to step $(6)$.}
\EndIf
\ElsIf {$t_{n-1}+k_{n}==T$} {stop the algorithm;} 
\Else \State{let $k_{n}\gets\kappa k_{n}$, and set $n\gets n+1$.}
\EndIf
\EndWhile
\end{algorithmic}
\end{algorithm}

\begin{remark}
As already emphasized in Remark~\ref{im}, step $(6)$ in the above Algorithm~\ref{spacetime} may be delicate to realize if derefinements of the spatial mesh are taken into account. Indeed, in step~(5), any coarsening procedure should be moderate in order to prevent the Newton iteration from leaving the current basin of attraction. 

We also remark that the parameter $ k_{\min}>0$ in step~(4) guarantees that the step size $k_{n}$ does not become overly small. This restriction needs to be relaxed when resolving \emph{finite time blow-up problems}, where the adaptivity with respect to the time evolution requires arbitrarily small step sizes~$k_n$ close to the blow-up time; see, e.g., \cite{JanssenWihler:15,BandleBrunner:98,Nakagawa:75} for details.
\end{remark}


\section{Numerical Experiments}\label{sec:numerics}
We will now illustrate and test the above Algorithm~\ref{spacetime} by means of a number of numerical experiments. We choose the initial spatial meshes to be sufficiently fine (as to fulfill~\eqref{eq:eta0}). Elements~$K\in \mathcal{T}_{h}^{n} $ are derefined whenever $ \eta_{h,K,N}<0.1 \overline{\eta}_{h,\Omega,N} $, where $ \overline{\eta}_{h,\Omega,N}$ signifies the mean of all~$\eta_{h,K,N}$, $K\in\mathcal{T}^n_h$; see~\eqref{eq:eta}.
 
\begin{example}\label{ex:1}
On~$\Omega=(0,1)$ let us consider the {\em linear} singularly perturbed initial/boundary value problem
\[
\begin{aligned}
\partial_{t}u -\varepsilon u''&=\exp(t) \ &&\text{on} \ \Omega\times(0,T],\\
u&=0 \ &&\text{on} \ \{0,1\}\times(0,T], \\
u(0,\cdot)&=g_{\varepsilon}&&\text{in}\ \Omega,
\end{aligned}
\]
where $g_{\varepsilon}$ is the solution of the elliptic boundary value problem
\[
-\varepsilon g_{\varepsilon}'' + g_{\varepsilon} = 1, \qquad
g_\varepsilon(0)=g_\varepsilon(1)=0. 
\]
Note that $g_{\varepsilon}$ exhibits boundary layers at $ x\in \{0,1\}$; cf.~\cite{AmreinWihler:15}. Since the problem is linear, the Newton iteration is redundant in this example. We prescribe the time derefinement/refinement parameters $ \kappa =2,  \sigma=\nicefrac{1}{2}$. Moreover, we compute a numerical solution up to the final time $T=1$, and set the local error tolerances from \eqref{eq:localtolerances} (as well as $ \varepsilon_{0} $ given in \eqref{eq:eta0}) to $ 10^{-3}$. Furthermore, the initial time step $k_0$ is chosen to be~$ \nicefrac{1}{10}$.

Our goal here is to test the robustness of the \emph{a posteriori} error analysis with respect to $\varepsilon$ as $\varepsilon\to 0$. To this end, we quantify the performance of our algorithm by comparing the true error $ \norm{u-u_{\mathcal{I}}}_{L^{2}(0,t_n;Z)}^{2}+\norm{u-u_{\mathcal{I}}}_{L^{\infty}(0,t_n;L^{2}(\Omega))}^{2} $ 
with the estimated error (i.e., the right-hand side of~\eqref{eq:En}), and compute the time-dependent efficiency indices (defined by the ratio of the estimated error and the true error $ \norm{u-u_{\mathcal{I}}}_{L^{2}(0,t_n;Z)}^{2}+\norm{u-u_{\mathcal{I}}}_{L^{\infty}(0,t_n;L^{2}(\Omega))}^{2} $ for $ n \in \{1,2,\ldots,M\}$); the results are displayed in Figure~\ref{bild2} for $ \varepsilon = 10^{-p}$, with $p \in \{1,2,3,4,5\} $. They show that the boundary layers close to~$\{0,1\}$ are properly resolved, and clearly highlight the robustness of the efficiency indices with respect to~$\varepsilon$ as $\varepsilon \to 0$. 

\begin{figure}
\includegraphics[width=0.47\textwidth]{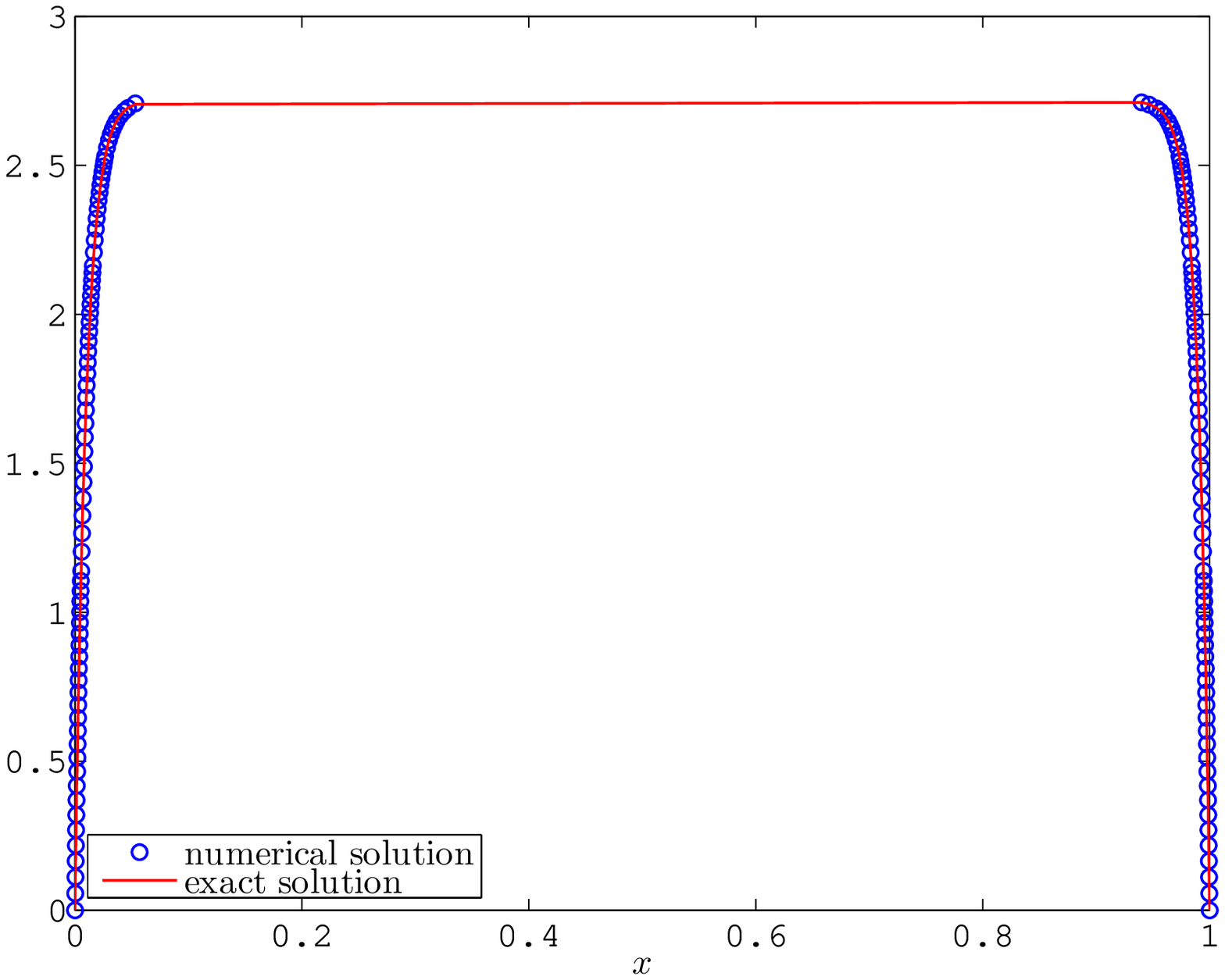}
\hfill
\includegraphics[width=0.46\textwidth]{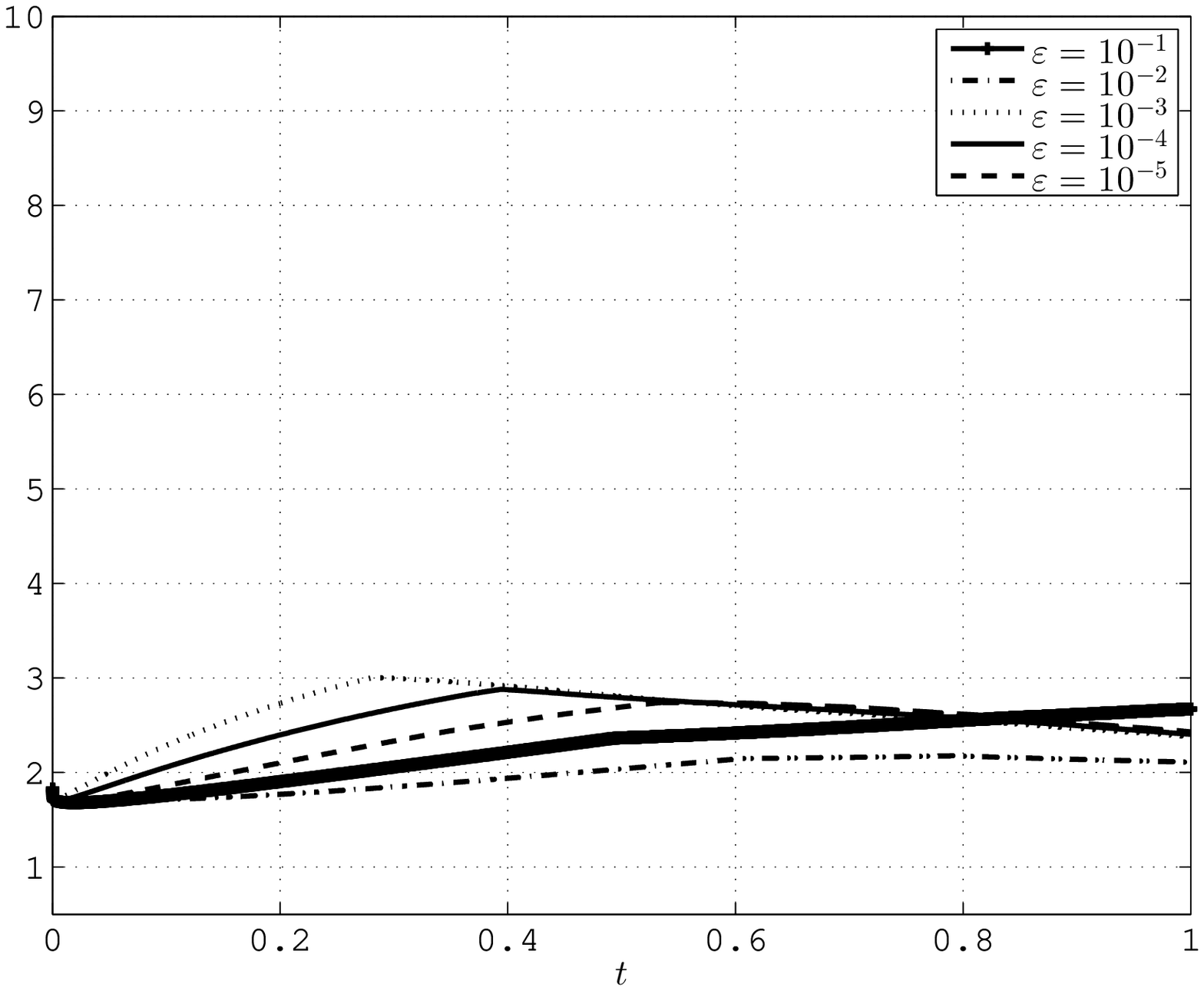}
\caption{Example~\ref{ex:1}: Numerical solution (with~'$\circ$' indicating the mesh points) vs. exact solution at $T=1$ (left), and efficiency indices for various choices of $\varepsilon$ (right).}
\label{bild2}
\end{figure}
\end{example}

\begin{example}\label{ex:2}
Furthermore, on~$\Omega = (0,1)$ consider the following nonlinear singularly perturbed initial/boundary value problem:
\begin{equation}
\label{boundarylayers}
\left\{ \begin{aligned}
\partial_{t}u - \varepsilon u'' &= -u^{4}+\sin(t),   && \text{on} \  \Omega\times(0,T],\\
u&=0,  &&\text{on} \  \{0,1\}\times(0,T],\\
u(0,\cdot)&=g  && \text{in}  \ \Omega.
\end{aligned} \right.
\end{equation}
When evolving in time, problem \eqref{boundarylayers} exhibits boundary layers for $ 0<\varepsilon \ll 1 $; see Figure ~\ref{bild22} (left), and \cite{Verhulst} for a detailed discussion of this problem. We consider~$\varepsilon=10^{-5}$, and choose the local error tolerances from \eqref{eq:localtolerances} (as well as $\varepsilon_{0}$ from~\eqref{eq:eta0}) to be~$10^{-3}$, and the initial time step as~$k_{0}=\nicefrac{1}{4}$. In Figure~\ref{bild22} (right) we depict a log/log plot of the estimated error from \eqref{eq:En} up to the final time $T=2$. Notice that the slope~$\nicefrac{1}{2}$ in the log/log plot is due to the fact that
\begin{align*}
\sqrt{E^{n}(u_{\mathcal{I}},g)} & \preccurlyeq \left(\eta_0^2+\sum_{l=1}^nk_l\left(\eta_{l,\Omega,N}^2+\vartheta^2_{l,\Omega,N}+\Upsilon^2_{l,\Omega,N}\right)\right)^{\nicefrac{1}{2}}\\
& \leq \left(\varepsilon_{0}^{2}+t_n(\varepsilon_{\text{loc},\eta}^{2}+\varepsilon_{\text{loc},\vartheta}^{2}+\varepsilon_{\text{loc},\Upsilon}^{2})\right)^{\nicefrac{1}{2}},
\end{align*}
i.e., for sufficiently small~$\varepsilon_0>0$, we expect the error to grow of order~$ \mathcal{O}(t_n^{\nicefrac12})$ as time evolves.

\begin{figure}
\includegraphics[width=0.479\textwidth]{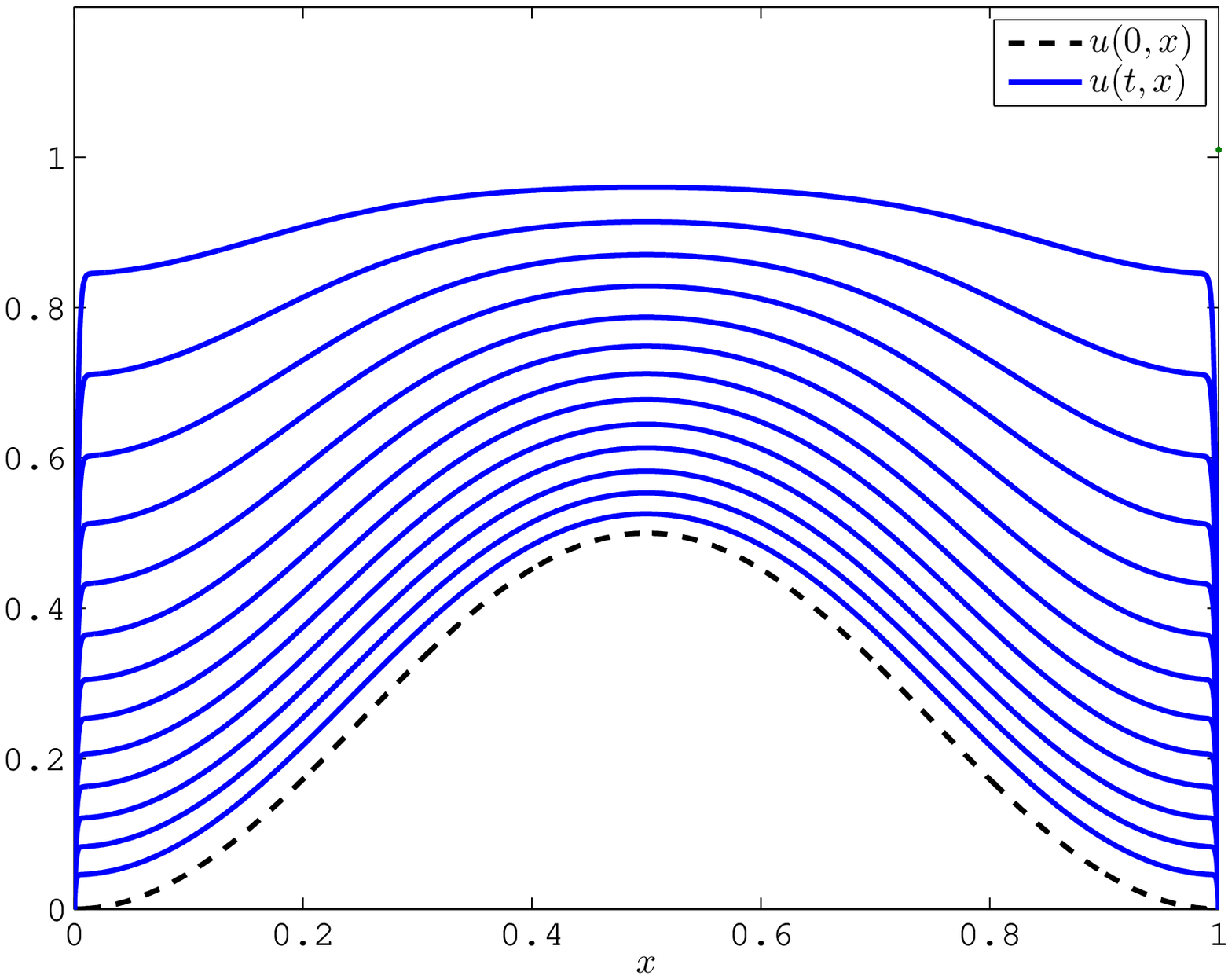}
\hfill
\includegraphics[width=0.47\textwidth]{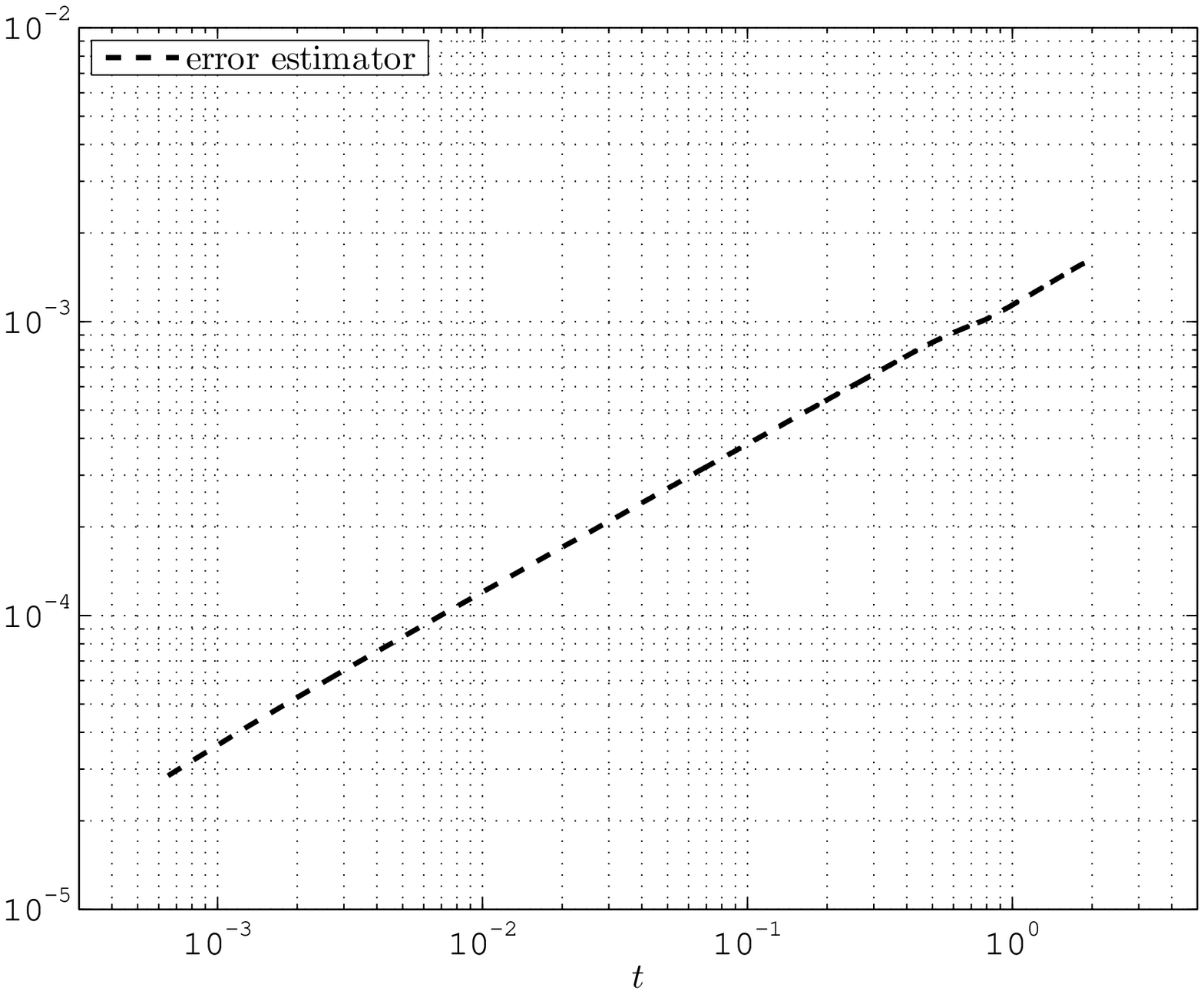}
\caption{Example~\ref{ex:2}: Snapshots of the numerical solution (as time is evolving) corresponding to problem \eqref{boundarylayers} with $\varepsilon =  10^{-5}$  (left), and the estimated error for $\sqrt{E^{n}(u_{\mathcal{I}},g)}$ (right).}
\label{bild22}
\end{figure}

\end{example}

\begin{example}\label{ex:3}
Finally, we consider the nonlinear problem
\begin{equation}
\label{blow up}
\left\{ \begin{aligned}
\partial_{t}u - \varepsilon \Delta u &= u^{\beta},   && \text{on} \ \Omega\times(0,T],\\
u&=0,  &&\text{on} \   \{0,1\}\times(0,T],\\
u(0,\cdot)&=g  && \text{in}  \  \Omega.
\end{aligned} \right.
\end{equation}
A detailed discussion of problems with power-type source terms can be found, for instance, in the monograph~\cite{SamarskiiGalaktionovKurdyumovMikhailov:95}. In particular, for~$\beta>1$, the solution of~\eqref{blow up} will become unbounded  in finite time provided the initial data $u(0,\cdot)=g\geq 0 $ is suitably chosen.

On the left of Figure~\ref{bild1} we show the numerical solution of~\eqref{blow up} for $ \Omega = (0,4)$, $ \beta=4$, $\varepsilon=10^{-3}$, and $T\approx 0.1$. The local error tolerances from \eqref{eq:localtolerances} are set to~$10^{-2}$, and $k_{0}=10^{-3}$. On the right in Figure~\ref{bild1}, we present a log/log plot of the estimated error (i.e. the right-hand side of \eqref{eq:En}) corresponding to the numerical solution shown on the left in Figure~\ref{bild1}. We clearly observe that the estimated error from \eqref{eq:En} increases with slope $=\nicefrac{1}{2}$, as in Example~\ref{ex:2}.
Moreover, as time evolves, we see that the adaptive procedure is able to resolve properly the spike located around $x=2$. 

\begin{figure}[!ht]
\includegraphics[width=0.474\textwidth]{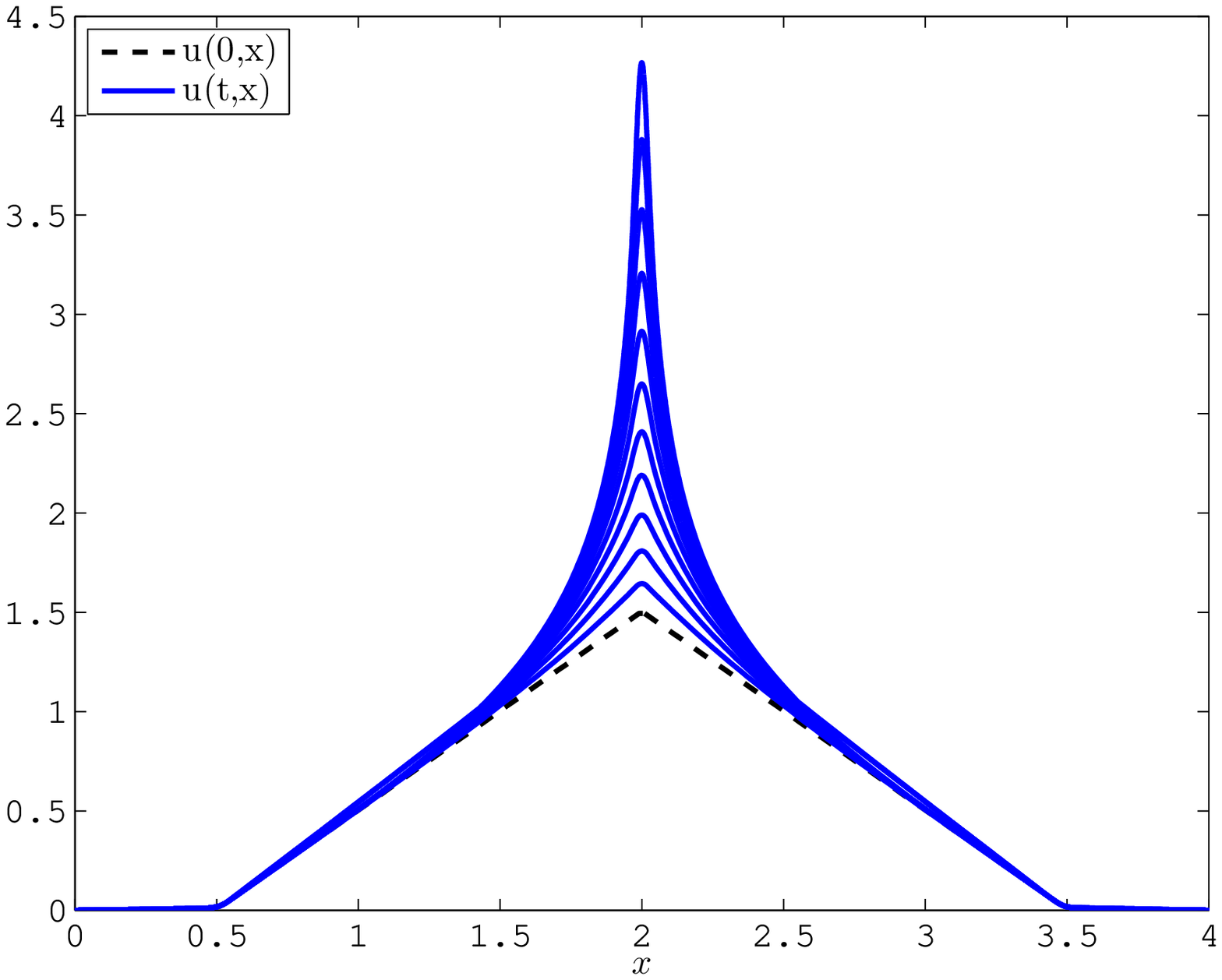}
\hfill
\includegraphics[width=0.471\textwidth]{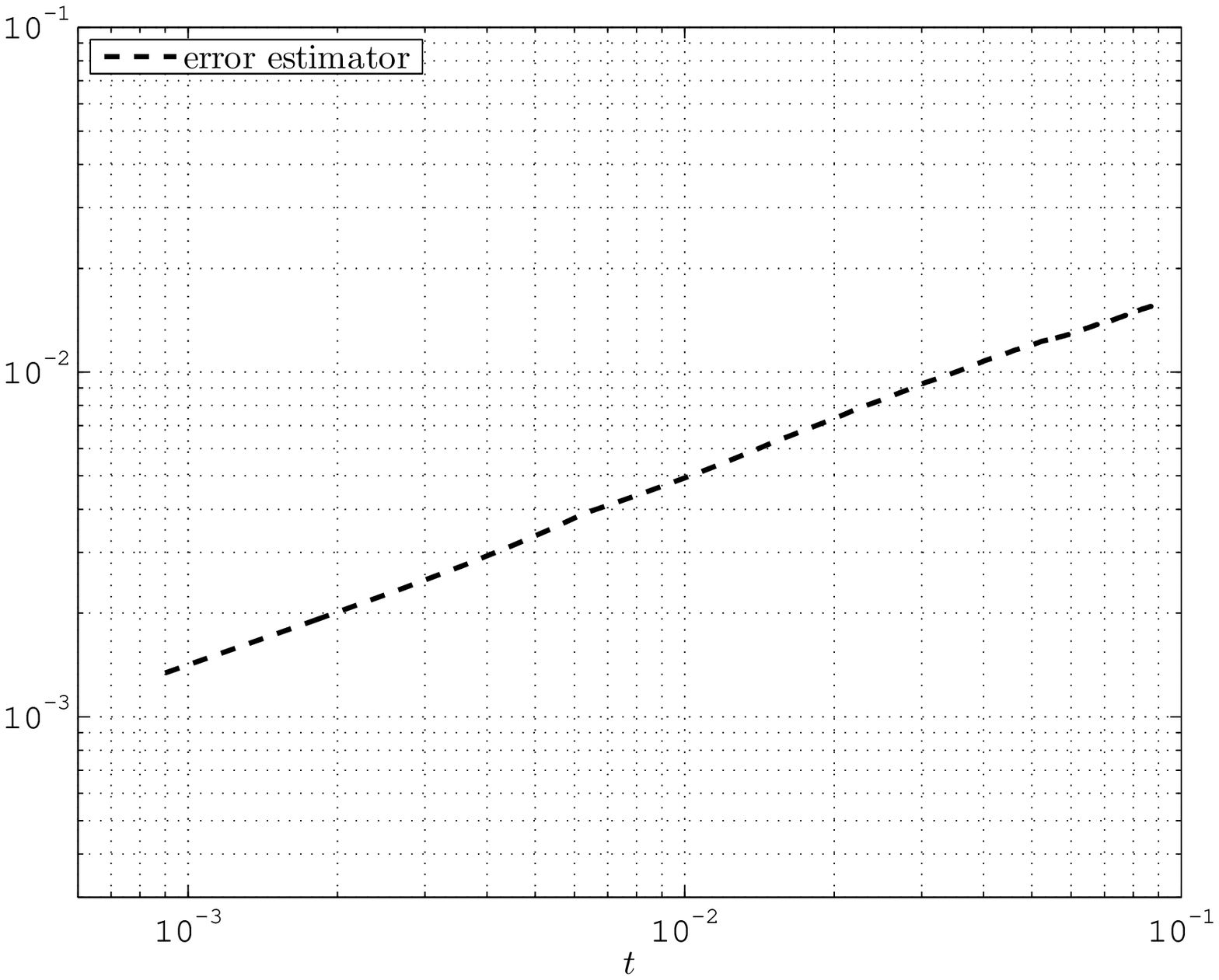}
\caption{Example~\ref{ex:3}: Snapshots of the numerical solution (as time is evolving) corresponding to problem \eqref{blow up} with $\varepsilon=10^{-3}
$, $\beta = 4 $, $ \Omega = [0,4] $ (left), and the estimated error for $\sqrt{E^{n}(u_{\mathcal{I}},g)}$ (right).}
\label{bild1}
\end{figure}
\end{example}

\section{Conclusions}\label{sc:conclusions}

The aim of this paper is the development of a reliable and computationally efficient procedure for the numerical solution of semilinear parabolic boundary value problems with possible singular perturbations. The key idea is to employ Newton's method to locally linearize the problem, and to apply an automatic (spatial) finite element mesh refinement approach as well as an adaptive time stepping control procedure. The numerical scheme is studied within the context of a robust 
(with respect to the singular perturbations) {\em a posteriori} 
residual-oriented error analysis, and a corresponding adaptive mesh refinement scheme is developed. Our numerical experiments clearly illustrate the ability of the proposed methodology to reliably find solutions, and to robustly resolve the singular perturbations at the expected rate.

\bibliographystyle{amsplain}
\bibliography{references}

\end{document}